\documentclass[10pt]{amsart}
\usepackage{latexsym, amsmath,amssymb}

\renewcommand{\epsilon}{\varepsilon}
\newcommand{\newsection}[1]
{\subsection{#1}\setcounter{theorem}{0} \setcounter{equation}{0}
\par\noindent}

\newtheorem{theorem}{Theorem}

\newtheorem{lemma}[theorem]{Lemma}
\newtheorem{corr}[theorem]{Corollary}

\newtheorem{proposition}[theorem]{Proposition}
\newtheorem{deff}[theorem]{Definition}

\newcommand{\bth}{\begin{theorem}}
\newcommand{\ble}{\begin{lemma}}
\newcommand{\bcor}{\begin{corr}}

\newcommand{\bdeff}{\begin{deff}}

\newcommand{\bprop}{\begin{proposition}}
\newcommand{\ele}{\end{lemma}}
\newcommand{\ecor}{\end{corr}}
\newcommand{\edeff}{\end{deff}}

\newcommand{\eprop}{\end{proposition}}

\newcommand{\cd}{\, \cdot\, }

\newcommand{\Rn}{{\mathbb R}^n}

\newcommand{\la}{\lambda}

\newcommand{\e}{\varepsilon}

\renewcommand{\Pi}{\varPi}
\renewcommand{\Re}{\rm{Re} \,}

\renewcommand{\epsilon}{\varepsilon}

\newcommand{\tidle}{\tilde}

\newcommand{\R}{{\mathbb R}}

\newcommand{\tube}{{{\mathcal T}_{\lambda^{-\frac12}}(\gamma)}}
\newcommand{\tubej}{{{\mathcal T}_{{\lambda_{j_k}}^{-\frac12}}(\gamma_{j_k})}}

\newcommand{\1}{{\rm 1\hspace*{-0.4ex}%
\rule{0.1ex}{1.52ex}\hspace*{0.2ex}}}

\newcommand{\vertiii}[1]{{\vert\kern-0.25ex\vert\kern-0.25ex\vert #1 
    \vert\kern-0.25ex\vert\kern-0.25ex\vert}}

\thanks{The author was supported in part by the NSF grant DMS-1361476.}

\begin{document}

\title[Problems related to the concentration of eigenfunctions]
{Problems related to the  \\ concentration of eigenfunctions}
%\thanks{The authors were supported in part by the NSF}
%
%
%
%
%
%
%\author{Jason Metcalfe}
\author{Christopher D. Sogge}
\address{Department of Mathematics,  Johns Hopkins University,
Baltimore, MD 21218}
\email{sogge@jhu.edu}

\begin{abstract}  We survey recent results related to the concentration of eigenfunctions.  We
also prove some new results concerning ball-concentration, as well as showing that eigenfunctions
saturating lower bounds for $L^1$-norms must also, in a measure theoretical sense, have
extreme concentration near a geodesic.
\end{abstract}
\maketitle

\newsection{Introduction}

In recent years much work has been devoted to problems related to the various types of
concentration exhibited by eigenfunctions of the Laplacian on Riemannian manifolds.
These include estimates for the their nodal sets, bounds for their restrictions to submanifolds,
especially geodesics, as well as estimates for their $L^2$-mass on small geodesic
tubes or balls.  In this paper we shall go over some of this work and also present a couple
of new results regarding the latter.  

Our eigenfunctions satisfy the equation
$$-\Delta_g e_\la(x)=\la^2 e_\la(x),$$
so that $\la$ is the frequency (and eigenvalue of $\sqrt{-\Delta_g}$), and we shall always
assume that they are $L^2$-normalized, i.e.,
$$\int_M|e_\la|^2 \, dV_g=1.$$
Here $\Delta_g$ and $dV_g$ are the Laplace-Beltrami operator and volume element, 
respectively, on our $n$-dimensional compact manifold $(M,g)$.

Let us first recall the extreme types of concentration that occur on the standard
round sphere, $S^n$.  On this manifold extreme concentration at points, as well as extreme
concentration along geodesics occur.  Both phenomena are present there due to the fact that
geodesic flow on $S^n$ is periodic and every geodesic is stable.

The eigenfunctions that have maximal concentration at points on $S^n$ are the zonal functions. 
Recall that the eigenvalues of $\sqrt{-\Delta_{S^n}}$ are $\sqrt{k(k+n-1)}$, repeating with 
multiplicity $d_k\approx k^{n-1}$.  If ${\mathcal H}_k$ is the associated space of spherical
harmonics of degree $k$ and $\{e_{k,1},\dots e_{k,d_k}\}$ an orthonormal basis of ${\mathcal H}_k$,
then the $L^2$-normalized zonal function centered about a point $x_0\in S^n$ is given by
the formula
$$Z_\la(x)=\Biggl(\frac{|S^n|}{d_k}\Biggr)^{\frac12}\sum_{j=1}^{d_k}e_{k,j}(x)\overline{e_{k,j}(x_0)},
\quad \la = \sqrt{k(k+n-1)},$$
with $|S^n|$ denoting the volume of $S^n$.  Usually one takes $x_0$ to be the north
pole $\1=(0,\dots,0,1)$ if we write $S^n$ as $\{x\in {\mathbb R}^{n+1}: \, \sum x_j^2=1\}$.
For this choice of $x_0$, by the classical Darboux-Szeg\"o formula, modulo lower order
terms, if $d_g(\cd, \cd)$ denotes geodesic distance, then
\begin{multline}\label{1.1}
Z_\la(x)\approx \frac{\cos\bigl[(k+\tfrac{n-1}2) d_g(x,\pm \1)+\tfrac{(n-1)\pi}4\bigr]}{(d_g(x,\pm \1))^{\frac{n-1}2}} \, ,
\quad
\text{if} \, \, \, d_g(x,\pm \1)\ge \la^{-1}, \\ 
\text{with} \, \, \, \la=\sqrt{k(k+n-1)},
\end{multline}
as well as
\begin{equation}\label{1.2}
Z_\la(x)=O(\la^{\frac{n-1}2}), \quad \text{if } \, \, d_g(x,\pm \1)\le \la^{-1}.
\end{equation}
Since $Z_\la(\1)=(d_k/|S^n|)^{\frac12}\approx \la^{\frac{n-1}2}$, we deduce from the above that
$\|Z_\la\|_\infty \approx \la^{\frac{n-1}2}$, and if we use \eqref{1.1}--\eqref{1.2}, we conclude
further that
\begin{equation}\label{1.3}
\|Z_\la\|_{L^p(S^n)}\approx \la^{n(\frac12-\frac1p)-\frac12}, \quad \text{if } \, \, p>\tfrac{2n}{n-1}.
\end{equation}

The other type of extreme concentration on $S^n$ is exhibited by the highest weight spherical
harmonics,
\begin{equation}\label{1.4}
Q_\la(x)=c_k \, \, {\Re }(x_1+ix_2)^k, \quad \la=\sqrt{k(k+n-1)},
\end{equation}
where
\begin{equation}\label{1.5}
c_k\approx k^{\frac{n-1}4},
\end{equation}
so that $\int |Q_\la|^2=1$.  Since
$$|x_1+ix_2|^k\approx e^{-\frac{k}2 |x'|^2}, \quad \text{on } \, \, S^n,
\quad \text{if } \, \, x'=(x_3,\dots,x_{n+1}),$$
we see that the mass of $Q_\la$ is concentrated on a $\la^{-\frac12}$ tubular
neighborhood of the equator, $\gamma$, where $x'=0$.  From the above one can
deduce that
\begin{equation}\label{1.6}
\|Q_\la\|_{L^p(S^n)}\approx \la^{\frac{n-1}2(\frac12-\frac1p)}, \quad p>2,
\end{equation}
as well as 
\begin{equation}\label{1.7}
\|Q_\la\|_{L^1(S^n)} \approx \la^{-\frac{n-1}4},
\end{equation}
and 
\begin{equation}\label{1.8}
\|Q_\la\|_{L^2(\tube)}\approx 1,
\end{equation}
if $\tube$ denotes a $\la^{-\frac12}$ tubular neighborhood of $\gamma$.

In \cite{SSh} it was shown that the $L^p$-bounds in \eqref{1.3} and \eqref{1.6} are extreme among
the spherical harmonics.  In \cite{Sef}, it was also shown that they are the worst possible
for $L^p$-norms on $n$-dimensional compact Riemannian manifolds.  Specifically, the 
$L^p$-upper bounds,
\begin{equation}\label{1.9}
\|e_\la\|_{L^p(M)}\lesssim \la^{\sigma(p)}, \quad 2<p\le \infty,
\end{equation}
were obtained with
\begin{equation}
\label{1.10}
\sigma(p)=\max\bigl\{(\tfrac{n-1}2(\tfrac12-\tfrac1p), \, n(\tfrac12-\tfrac1p)-\tfrac12\bigr\}.
\end{equation}
Thus, since $\tfrac{n-1}2(\tfrac12-\tfrac1p)=n(\tfrac12-\tfrac1p)-\tfrac12$ when $p=\tfrac{2(n+1)}{n-1}$,
the bounds in \eqref{1.9} can be rewritten as
\begin{equation}\label{1.9'}\tag{1.9$'$}
\|e_\la\|_{L^p(M)}\lesssim
\begin{cases} \la^{\frac{n-1}2\bigl(\frac12-\frac1p\bigr)}, \quad 2\le p\le \frac{2(n+1)}{n-1},
\\ \\
\la^{n(\frac12-\frac1p)-\frac12}, \quad \frac{2(n+1)}{n-1} \le p\le \infty.
\end{cases}
\end{equation}
More, generally, it was shown in \cite{Sef} that one has the following estimates
for the unit-band spectral projection operators
\begin{equation}\label{1.9''}\tag{1.9$''$}
\|\chi_{[\la,\la+1]}\|_{L^2(M)\to L^p(M)}\lesssim 
\begin{cases} \la^{\frac{n-1}2\bigl(\frac12-\frac1p\bigr)}, \quad 2\le p\le \frac{2(n+1)}{n-1},
\\ \\
\la^{n(\frac12-\frac1p)-\frac12}, \quad \frac{2(n+1)}{n-1} \le p\le \infty.
\end{cases}
\end{equation}
If $I\subset \R_+$ is an interval, then $\chi_I$ denotes the projection onto
frequencies in $I$.  Thus, if $0=\la_0<\la_1\le \la_2\le \dots$ are the eigenvalues
of $\sqrt{-\Delta_g}$ counted with multiplicity and if $\{e_{\la_j}\}$ is the associated orthonormal
basis of eigenfunctions, then
\begin{equation*}
\chi_If=\sum_{\la_j\in I}E_jf, \quad \text{with } \, \,
E_jf(x)=e_{\la_j}(x) \times \int_M f(y)\, \overline{e_{\la_j}(y)}\, dV_g.
\end{equation*}

The estimate \eqref{1.9''} is a ``local estimate" that cannot be improved on any
Riemannian manifold (see \cite{SFIO}).  However, in many cases, the estimate \eqref{1.9'}
can be improved, even though it is saturated on $S^n$.

\newsection{Improved $L^p$ bounds for large exponents}

In 2002, the author and Zelditch~\cite{SZ} showed that for generic Riemannian manifolds
the estimates in \eqref{1.9} can be improved when $p$ is larger than the ``critical" index
$p_c=\tfrac{2(n+1)}{n-1}$.  Specifically, if for a given $x\in M$, and an initial unit direction
$\xi\in S_x^*M$ over $x$, we say that $\xi\in {\mathcal L}_x$ if the geodesic with initial
direction $\xi$ loops back through $x$ in some positive time $t$, then the following
result was obtained in \cite{SZ}:

\begin{theorem}  If $|{\mathcal L}_x|=0$ for all $x\in M$, then
$\|e_\la\|_{L^\infty(M)}=o(\la^{\frac{n-1}2})$, and hence if $\sigma(p)$ is as in
\eqref{1.10}
\begin{equation}\label{2.1}
\|e_\la\|_{L^p(M)}=o(\la^{\sigma(p)}), \quad \forall \, p>p_c=\tfrac{2(n+1)}{n-1}.
\end{equation}
\end{theorem}

The bounds in \eqref{2.1} follow via interpolation from the estimate $\|e_\la\|_{L^\infty(M)}=o(\la^{\frac{n-1}2})$ and the special case of \eqref{1.9'} with $p=p_c$.  It was shown in
\cite{SZ} that the condition $|{\mathcal L}_x|=0$ for all $x\in M$ holds for generic $(M,g)$.
Here $|\cd|$ denotes the measure on $S^*_xM$ coming from restricting Liouville measure
on $T^*M$ to $S^*_xM$.

The above theorem was improved by the author, Toth and Zelditch~\cite{STZ} by showing that
the condition that the looping directions is of measure zero for all $x$ can be replaced
by the condition that the set of recurrent directions, ${\mathcal R}_x\subset S^*_xM$,
be of measure zero for every $x\in M$.

Further improvements were obtained recently by the author and Zelditch~\cite{SZRA}--\cite{SZRA2}
under the assumption that $(M,g)$ is real analytic.  Indeed, in this case, necessary and
sufficient conditions for improved sup-norm estimates were obtained. 

 In two dimensions,
the result is very natural and simple to state:

\begin{theorem}  Let $(M,g)$ be a real analytic two-dimensional Riemannian surface.  Then
\begin{equation}\label{2.2}
\|e_\la\|_{L^\infty(M)}=o(\la^{\frac12})
\end{equation}
if and only if there is no point $x$ through which the geodesic flow is periodic.  In this case,
\eqref{2.1} also holds.
\end{theorem}

There is also a result in higher dimensions that is a bit more technical to state.   Based on an
earlier work \cite{SZ}, it was known that, in the real analytic case, if $|{\mathcal L}_x|\ne 0$ for
some $x$, then one must have that ${\mathcal L}_x=S^*_xM$.  In other words, $x$ must
be a ``focal point,'' and furthermore, in this case, there must be a minimal positive  $T_0$ so that
all unit speed geodesics loop back in time $T_0$.  As a result, unitary Perron-Frobenius operators
on $L^2(S^*_xM)$ are well-defined, and the theorem says that one has \eqref{2.2}
if and only if there are no focal points for which the associated Perron-Frobenius operator has
a non-trivial invariant function.  These operators were introduced earlier to study these types of problems
by Safarov~\cite{Sa}, and our analysis uses ideas from his work.

If one makes curvature assumptions, then one can get further improvements.  Indeed, it has
been known for some time that if the sectional curvatures of $(M,g)$ are nonpositive
then one has
\begin{equation*}
\|e_\la\|_{L^\infty(M)}=O\bigl(\la^{\frac{n-1}2}/\sqrt{\log \la}\bigr).
\end{equation*}
This follows directly from B\'erard's \cite{PB} proof of improved error term estimates
for the Weyl formula under these curvature assumptions (see, e.g., \cite{SoHang}).  Recently,
Hassell and Tacy~\cite{HT} were able to extend this by showing that one has
\begin{equation*}
\|e_\la\|_{L^\infty(M)}=O\bigl(\la^{\sigma(p)}/\sqrt{\log \la}\bigr),
\quad \text{if } \, \, p>p_c=\tfrac{2(n+1)}{n-1}.
\end{equation*}

A very interesting open problem would be to show that one also has these sorts of improvements
for the endpoint $p=p_c$ for all eigenfunctions.  Recently, under the assumption of negative
sectional curvatures, Hezari and Rivi\`ere~\cite{HR} were able to obtain improvements
involving different powers of $\log\la$, but for eigenfunctions associated with a density
one subsequence of eigenvalues.

\newsection{Kakeya-Nikodym norms and
improved $L^p$ bounds for small exponents}

As we mentioned before, the highest weight spherical harmonics, $Q_\la$, on $S^n$ saturate $L^p$-norms
for $2<p<p_c=\tfrac{2(n+1)}{n-1}$.  They have significant $L^2$ mass on tubes $\tube$ of width $\la^{-\frac12}$ 
about unit geodesics $\gamma$ on the equator.  Given a compact Riemannian manifold $(M,g)$ we shall let
$\varPi$ denote the space of unit length geodesics and $\tube$ the resulting $\lambda^{-\frac12}$-tube about
$\gamma$.  Then the highest weight spherical harmonics also saturate the Kakeya-Nikodym norms given by
\begin{equation}\label{3.1}
{\vertiii{e_\la}}_{KN}=\sup_{\gamma\in \varPi} \|e_\la\|_{{L^2(\tube)}},
\end{equation}
in the sense that
\begin{equation}\label{3.2}
\liminf_{\la\to \infty}\vertiii{Q_\la}_{KN}>0.
\end{equation}
Note that, since we are assuming that our eigenfunctions have $L^2(M)$-norm one, we always have the trivial upper bound
\begin{equation}\label{3.3}
\vertiii{e_\la}_{KN}\le 1.
\end{equation}

These norms were introduced by the author~\cite{SKN} as a way to control $L^p$-norms in two-dimensions
for exponents $2<p<p_c=6$, following earlier related work of Bourgain~\cite{Bo}.  

The author was able to obtain the following result.

\begin{theorem}  Let $(M,g)$ be a two-dimensional compact Riemannian manifold.  Then if $\{\la_{j_k}\}$ is a subsequence
of eigenvalues, the following are equivalent
\begin{gather}
\la_{j_k}^{-\sigma(p)}\|e_{\la_{j_k}}\|_{L^p(M)}\to 0,\quad \forall \, 2<p<6 \label{3.4}
\\
\vertiii{e_{\la_{j_k}}}_{KN}\to 0 \label{3.5}
\\
\la_{j_k}^{-\frac14}\sup_{\gamma\in \varPi}\Bigl(\int_\gamma |e_{\la_{j_k}}|^2 \, ds\Bigr)^{\frac12}\to 0 \label{3.6}.
\end{gather}
\end{theorem}

Burq, G\'erard and Tzvetkov~\cite{BGT} showed that in two-dimensions one has
\begin{equation}\label{3.7}
\sup_{\gamma\in \varPi}\Bigl(\int_\gamma |e_\la|^2 \, ds\Bigr)^{\frac12} \lesssim \la^{\frac14}\|e_\la\|_{L^2(M)},
\end{equation}
which is another estimate saturated by the $Q_\la$.  Bourgain~\cite{Bo} showed that, more generally, one has
\begin{equation}\label{3.8}
\sup_{\gamma\in \varPi} \Bigl(\int_\gamma |e_\la|^2 \, ds\Bigr)^{\frac12} \lesssim \la^{\frac1{2p}}\|e_\la\|_{L^p(M)},
\quad 2\le p\le \infty,
\end{equation}
which gives that \eqref{3.4} implies \eqref{3.6}.  Clearly, \eqref{3.6} implies \eqref{3.5}.  The fact that 
\eqref{3.5} implies \eqref{3.4} was established in \cite{SKN} where inequalities were obtained which imply that
\begin{equation}\label{3.9}
\|e_\la\|_{L^4(M)}\lesssim \vertiii{e_\la}_{KN}^{\frac14}.
\end{equation}
This of course implies that whenever \eqref{3.5} is valid we must have the special case of \eqref{3.4}
corresponding to $p=4$, and the corresponding results for the other exponents $2<p<6$ follow via
interpolation with the estimate \eqref{1.9'} for $p=6$ and the trivial bound for $p=2$.

In \cite{SKN}, using microlocal analysis, it was shown that on any two-dimensional manifold one has
$$\Bigl(\int_\gamma |e_\la|^2 \, ds \Bigr)^{\frac12}=o(\la^{\frac14})$$
provided that $\gamma\in \varPi$ is not a unit-segment of a periodic geodesic.

The author and Zelditch used the Hadamard parametrix and microlocal analysis in \cite{SZKN} to show
that if $(M,g)$ is a two-dimensional manifold of nonpositive curvature then \eqref{3.6} is valid for
all eigenfunctions, i.e.,
\begin{equation}\label{3.10}
\sup_{\gamma\in \varPi}\Bigl(\int_\gamma |e_\la|^2 \, ds \Bigr)^{\frac12}=o(\la^{\frac14}).
\end{equation} 
Improvements were obtained by Chen and the author~\cite{CS} in three dimensions as well under the assumption of constant
nonpositive sectional curvature.

In \cite{BS2} Blair and the author found the appropriate extensions of these results to higher dimensions.  It turns
out that, for $n\ge4$, $L^2(\gamma)$ restriction estimates are too singular to control $L^p(M)$ norms for $2<p<p_c=\tfrac{2(n+1)}{n-1}$.
This is due in part to the fact that, for these dimensions, the zonal functions, $Z_\la$, saturate $L^2(\gamma)$ norms,
but not these $L^p(M)$ norms, while the highest weight spherical harmonics, $Q_\la$, saturate the latter but not the former.
Notwithstanding, the less singular Kakeya-Nikodym norms control $L^p(M)$ norms for $2<p<p_c$, and it was shown in
\cite{BS2} that 
$$\la_{j_k}^{-\sigma(p)}\|e_{\la_{j_k}}\|_{L^p(M)}\to 0,\quad \forall \, \, 2<p<p_c=\tfrac{2(n+1)}{n-1},$$
if and only if
$$\vertiii{e_{\la_{j_k}}}_{KN}\to 0.$$
Additionally, it was shown in \cite{BS2} that, in any dimension, if $(M,g)$ is a manifold with nonpositive sectional curvatures, then 
\begin{equation}\label{3.10'}\tag{3.10$'$}
\vertiii{e_\la}_{KN}=o(1),
\end{equation}
and hence 
\begin{equation}\label{3.11}
\|e_\la\|_{L^p(M)}=o(\la^{\sigma(p)}), \quad 2<p<\tfrac{2(n+1)}{n-1}.
\end{equation}

These results have been refined considerably in recent papers of Blair and the author~\cite{BS1}, \cite{BS2} and \cite{BSlog}.
The linkage between $L^p(M)$ norms and Kakeya-Nikodym estimates has been improved in \cite{BS15}, where the following
estimate was established
\begin{equation}\label{3.12}
\|e_\la\|_{L^p(M)}\lesssim \la^{\frac{n-1}2(\frac12-\frac1p)}\vertiii{e_\la}_{KN}^{\frac{2(n+1)}{n-1}(\frac1p-\frac{n-1}{2(n+1)})},
\quad \text{if } \, \, \tfrac{2(n+2)}n<p<\tfrac{2(n+1)}{n-1}.
\end{equation}

Additionally, in \cite{BSlog} the following result was obtained
\begin{theorem}\label{mainthm}  Suppose $(M,g)$ has nonpositive sectional curvatures.  Then
\begin{equation}\label{3.13}
\sup_{\gamma\in \varPi}\int_\tube |e_\la|^2 \, dV \lesssim c(\lambda),
\end{equation}
for $\la\gg 1$ with
$$c(\lambda)=
\begin{cases}
(\log \la)^{-\frac12}, \quad \text{if  } \, n=2
\\
(\log \la)^{-1} \log\log \la, \quad \text{if  }\, n=3
\\
(\log \la)^{-1}, \quad \text{if  }\, n\ge 4.
\end{cases}
$$
Moreover, if $n=2$, we have
\begin{equation}\label{3.14}
\sup_{\gamma\in \varPi} \int_\gamma |e_\la|^2 \, ds \le C\la^{\frac12} c(\lambda).
%\le C(\la/\log \la)^{\frac12},
\end{equation}
\end{theorem}

By combining this result with \eqref{3.12}, one obtains the following
\begin{corr}\label{Lpcorr}  Assume, as above, that $(M,g)$ is a compact $n\ge2$ dimensional
manifold with nonpositive sectional curvatures.  Then for any $2<p<\tfrac{2(n+1)}{n-1}$
there is a number $\mu(p,n)>0$ so that
\begin{equation}\label{3.15}
\|e_\la\|_{L^p(M)}\lesssim \la^{\frac{n-1}2(\frac12-\frac1p)} \, \bigl(\log \la\bigr)^{-\mu(p,n)}.
\end{equation}
Furthermore, if $\tfrac{2(n+2)}n <p<\tfrac{2(n+1)}{n-1}$, one can take
\begin{multline}\label{3.16}   \quad
\mu(p,n)=
\begin{cases}
\tfrac{n+1}{n-1}(\tfrac1p-\tfrac{n-1}{2(n+1)}), \, \, \, \text{if } \, \,  n\ge 4,
\\ \\
\tfrac32(\tfrac1p-\tfrac16), \, \, \, \text{if } \, \, n=2, \end{cases}
\\ \text{and any } \, \, \mu(p,3)<2(\tfrac1p-\tfrac14), 
\, \, \text{if } \, \, n=3.
\end{multline}
\end{corr}

The proof of \eqref{mainthm} is similar to the argument from \cite{SZKN}.  A new ingredient, which allows the microlocal
analysis to obtain these logarithmic improvements under the assumption of nonpositive sectional curvatures, is the use
of the classical Topogonov triangle comparison theorem in Riemannian geometry.

\newsection{Lower bounds for $L^1$-norms and the size of nodal sets}

In \cite{SZnod1} the author and Zelditch showed that
there is a positive constant $c_M$ so that
\begin{equation}\label{4.1}
c_M\la^{-\frac{n-1}4}\le \|e_\la\|_{L^1(M)}.
\end{equation}
This is another estimate saturated by the highest weight spherical harmonics, $Q_\la$.  Hezari and the author~\cite{HS}
showed that if
$$Z_\la=\bigl\{x\in M: \, e_\la(x)=0\bigr\}$$
is the nodal set of a real eigenfunction, then its $(n-1)$-dimensional Hausdorff measure, $|Z_\la|$, satisfies
\begin{equation}\label{4.2}
\la \|e_\la\|_{L^1(M)}^2 \lesssim |Z_\la|.
\end{equation}
Combining \eqref{4.1} and \eqref{4.2} yields the lower bounds of Colding and Minicozzi~\cite{CM},
\begin{equation}\label{4.3}
\la^{1-\frac{n-1}2}\lesssim |Z_\la|,
\end{equation}
which is the best known lower bound for general smooth Riemannian manifolds.  Later,
the author and Zelditch in \cite{SZnod2} found a very simple proof of \eqref{4.3} by using
the estimate
\begin{equation}\label{4.1'}\tag{4.1$'$}
\|e_\la\|_{L^\infty(M)}\lesssim \la^{\frac{n-1}2}\|e_\la\|_{L^1(M)},
\end{equation}
which in their earlier paper, \cite{SZnod1}, was used to prove \eqref{4.1}. 

One also can use Kakeya-Nikodym norms to control lower bounds for $L^1(M)$ norms.  First, note that by
H\"older's inequality, if $p>2$,
$$1=\|e_\la\|_{L^2(M)}\le \|e_\la\|_{L^1(M)}^{\frac{p-2}{2(p-1)}} \,  \|e_\la\|_{L^p(M)}^{\frac{p}{2(p-1)}},$$
and so
$$\|e_\la\|_{L^1(M)}\ge \|e_\la\|_{L^p(M)}^{-\frac{p}{p-2}}=\la^{-\frac{n-1}4} \, \bigl(\la^{-\frac{n-1}2(\frac12-\frac1p)}\|e_\la\|_{L^p(M)}\bigr)^{-\frac{p}{p-2}}.
$$
Therefore, by \eqref{3.12},
\begin{equation}\label{4.4}
\la^{\frac{n-1}4}\|e_\la\|_{L^1(M)} \ge c_m\vertiii{e_\la}_{KN}^{-\frac{2(n+1)-(n-1)p}{(n-1)(p-2)}}, \quad
\text{if } \, \, \tfrac{2(n+2)}n <p<\tfrac{2(n+1)}{n-1}.
\end{equation}

Using this, we can deduce a couple of results.  The first, which is from \cite{BSlog}, is the following:

\begin{theorem}
Assume that $(M,g)$ is an $n$-dimensional compact manifold with nonpositive sectional curvatures.  Then
\begin{equation}\label{4.5}
\la^{-\frac{n-1}4} (\log\la)^{\mu} \lesssim \|e_\la\|_{L^1(M)},
\end{equation}
for any $\mu<\mu_n$ with
\begin{equation*}   \quad
\mu_n=
\begin{cases}
\tfrac{(n+1)^2}{n-1}, \quad \text{if } \, \, n\ge 3
\\ \\
\tfrac{(n+1)^2}{2(n-1)}, \quad \text{if } \, \, n=2.
\end{cases}
\end{equation*}
Consequently, if $e_\la$ is a real-valued eigenfunction and $|Z_\la|$ denotes
the $(n-1)$-dimensional Hausdorff measure of its nodal set, $Z_\la = \{x: \, e_\la(x)=0\}$, we have
\begin{equation}\label{4.6}
\la^{1-\tfrac{n-1}2} \bigl(\log \la\bigr)^{2\mu} \lesssim |Z_\la|,
\end{equation}
when $\mu<\mu_n$.  In particular, when $n=3$, $(\log \la)^{r}\lesssim |Z_\la|$ for
all $r<16$.
\end{theorem}

To prove this, one just uses the upper bound for $\vertiii{e_\la}_{KN}$ in \eqref{3.13}.  The resulting bounds
are optimized as $p\searrow\tfrac{2(n+2)}n$, which leads to \eqref{4.5}, which by \eqref{4.2}, implies \eqref{4.6}.

We can also use \eqref{4.4} to show that eigenfunctions saturating the lower bound \eqref{4.1} for $L^1(M)$-norms
must match the profile of the $Q_\la$ closely in the sense that there must be a geodesic tube $\tube$ on which a 
nontrivial portion of eigenfunction has values of size $\approx \la^{\frac{n-1}4}$.  Specifically, we have the following:

\begin{theorem}\label{scar}
Let $(M,g)$ be an $n$-dimensional compact Riemannian manifold.
Assume that $\la_{j_k}$ is a sequence  of eigenvalues for which
\begin{equation}\label{s.1}
\|e_{\la_{j_k}}\|_{L^1(M)}\le c_0\la_{j_k}^{-\frac{n-1}4}
\end{equation}
for some positive constant $c_0$.  Then there must be a positive
constant $c_1$ and tubes $\tubej$ centered
along a unit-length geodesic $\gamma_{j_k}\in \varPi$ 
%of some {\em periodic geodesic}
so that for suffiently large $\la_{j_k}$ we have
\begin{equation}\label{s.2}
c_1\le \int_{\tubej}|e_{\la_{j_k}}|^2 \, dV_g.
\end{equation}
Moreover, in this case, there also must be $\delta, \, c_2>0$,
depending only on $c_0$ and $(M,g)$, so
that for large $\la_{j_k}$ 
\begin{equation}\label{s.3}
{\rm Vol }\Bigl(\{x\in \tubej: 
|e_{\la_{j_k}}(x)|\in [c_2\la_{j_k}^{\frac{n-1}4}, \, c^{-1}_2\la_{j_k}^{\frac{n-1}4}]
\}\Bigr) \ge \delta \la_{j_k}^{-\frac{n-1}2}.
\end{equation}
\end{theorem}

To prove this, we note that \eqref{s.2} is an immediate consequence of \eqref{4.4}.
To prove the last part of the theorem we shall use \eqref{s.1} and
\eqref{4.1'}.
%the following inequality which was proved by the author and Zelditch~\cite{SZnod1} 
%\begin{equation*}
%\|e_\la\|_{L^\infty(M)}\lesssim \la^{\frac{n-1}2}\|e_\la\|_{L^1(M)}.
%\end{equation*}
We deduce, that if \eqref{s.1} is valid,  then we must have that there is
a constant $c_2$ so that
\begin{equation}\label{s.5}
\|e_{\la_{j_k}}\|_{L^\infty(M)}\le c_3\la_{j_k}^{\frac{n-1}4}.
\end{equation}
Thus, if we choose $c_2>0$ small enough so that $c_2^{-1}\ge c_3$ we 
have that $|e_{\la_{j_k}}(x)|\le c_2^{-1}\la_{j_k}^{\frac{n-1}4}$ on {\em all of}
$\tubej$.

To get the other condition in the left side of \eqref{s.3} on a nontrival
portion of $\tubej$ we need to use \eqref{s.2}.  We first note that
since $\text{Vol }(\tubej)\approx \la_{j_k}^{-\frac{n-1}2}$, there must
be a fixed constant $C_0$ so that
\begin{equation*}
\int_{ \bigl\{ x\in \tubej: \, |e_{\la_{j_k} }(x)|\le c_2\la_{j_k}^{\frac{n-1}4} \bigr\}}
|e_{\la_j}(x) |^2 \, dV_g \le C_0c_2.
\end{equation*}
If $c_1$ is as in \eqref{s.2}, we shall also assume that $c_2$ is chosen
small enough so that $C_0c_2\le c_1/2$.  We then have, by \eqref{s.5}
and \eqref{s.2}
\begin{multline*}
\frac{c_1}2\le \int_{\bigl\{x\in \tubej: \,
|e_{\la_j}(x)|\in [c_2\la_j^{\frac{n-1}4}, \, c^{-1}_2\la_j^{\frac{n-1}4}]
\bigr\}}|e_{\la_j}(x) |^2 \, dV_g
\\
\le c_2^{-2}\la_j^{\frac{n-1}2}
\text{Vol }\Bigl(
\{x\in \tubej: \,
|e_{\la_{j_k}}(x)|\in [c_2\la_j^{\frac{n-1}4}, \, c^{-1}_2\la_{j_k}^{\frac{n-1}4}]
\}
\Bigr).
\end{multline*}
From this we clearly get \eqref{s.3} if we take $\delta>0$ to be
$c_1c_2^2/2$, which completes the proof of Theorem \ref{scar}.

\newsection{ $L^2$-estimates for small balls}

In Han~\cite{H} and Hezari and Rivi\'ere~\cite{HR}, 
following in part earlier work of Zeldtich~\cite{ZelditchCMP},  it was shown that if $(M,g)$ has negative
sectional curvatures, then
\begin{equation}\label{5.1}
\int_{B(x,r)}|e_{\la_{j_k}}|^2 \, dV_g \, \le C_M \, r^n,\quad \text{if } \, \, r=\bigl(\log \la_{j_k}\bigr)^{\delta_n},
\end{equation}
for some $\delta_n>0$ depending on the dimension, as $\{\la_{j_k}\}$ ranges over a density one sequence
of eigenvalues so that the resulting system is quantum ergodic.\footnote{In Hezari and Rivi\'ere \cite[\S 3.1]{HR} uniform bounds of this type
were obtained for all balls $B(x_\ell,r)$, $\ell =1,\dots, N(r)$, $N(r)\approx r^{-n}$,  of  radius  $r$ as in \eqref{5.1} occurring in a specific covering of $M$.  The
covering can be chosen so that the  doubles, $B(x_\ell,2r)$, of the balls in their covering have uniformly
bounded overlap by an argument of \cite{CM}, which yields \eqref{5.1}.}  
Here $B(x,r)$ denotes the geodesic ball of
radius $r>0$ about $x\in M$.  Recent results of this type for toral eigenfunctions and much smaller balls
are due to Lester and Rudnick~\cite{LR}.

If $r$ does not depend on the frequency but is fixed, then a more precise version of \eqref{5.1} (also for a density
one sequence of eigenvalues) is a consequence of the Quantum Ergodic Theorem of Shnirelman~\cite{Snirelman} / Colin de Verdi\'ere~\cite{Colin} / Zelditch~\cite{Zelditch} under the assumption that the geodesic flow is ergodic (which is weaker than the assumptions in \cite{H} and \cite{HR}).  However,
obtaining results where the scale depends on the frequency is highly nontrivial.

Hezari and Rivi\'ere~\cite{HR} used \eqref{5.1} and a localization argument to show that under the assumption of negative curvature, for a density
one sequence of eigenfunctions one gets certain log-power improvements over the critical estimate \eqref{1.9'} with
$p=p_c=\tfrac{2(n+1)}{n-1}$.

The localization step was modified slightly by the author in \cite{Sball} by showing that for all $(M,g)$ (no curvature assumptions) one has
\begin{equation}\label{5.2}
\|e_\la\|_{L^{p_c}(M)}\le C\la^{\sigma(p_c)}\, \Bigl(r^{-\frac{n+1}4}\sup_{x\in M}\|e_\la\|_{L^2(B(x,r))}\Bigr)^{\frac2{n+1}},
\, \, \la^{-1}\le r\le {\rm Inj } \, M,
\end{equation}
with $p_c=\tfrac{2(n+1)}{n-1}$ and $\sigma(p)$ as in \eqref{1.10}.  Using this estimate one immediately sees how an estimate
like \eqref{5.1} leads to an improvement over \eqref{1.9'}.  \eqref{5.2} also shows, via the aforementioned Quantum Ergodicity Theorem, that if
the geodesic flow is ergodic then there must be a density one sequence of eigenfunctions for which
\begin{equation*}
\|e_{\la_{j_k}}\|_{L^{p_c}(M)}=o\bigl(\la_{j_k}^{\sigma(p_c)}\bigr).
\end{equation*}
A very intersting problem is whether this holds for all eigenfunctions assuming, say, negative or nonpositive sectional curvatures.

In \cite{Sball} the author showed that for any $(M,g)$ one always has
\begin{equation}\label{5.4}
\|e_\la\|_{L^2(B(x,r))}\le Cr^{\frac12}, \quad \la^{-1}\le r\le {\rm Inj } \, M,
\end{equation}
with $C$ depending only on $(M,g)$.  
By \eqref{1.1}--\eqref{1.2},
this estimate is saturated by the zonal functions, $Z_\la$, on $S^n$ for this full range of
$r$, and it is also saturated by the highest weight spherical harmonics, $Q_\la$, for the range $\la^{-1}\le r\le \la^{-\frac12}$.

Let us conclude this note by showing that for a slightly smaller range one can get log improvements over \eqref{5.4} under the assumption of 
nonpositive sectional curvatures:

\begin{theorem}\label{ball}
Fix $(M,g)$ with nonpositive sectional curvatures.  Then given $\e>0$ there is a constant $C_\e$ such that
\begin{equation}\label{5.5} 
\|e_\la\|_{L^2(B(x,r))}\le C_\e \bigl(r/\log \la)^{\frac12}, \quad \text{if } \, \, \la^{-1}\le r\le \la^{-\frac12-\e}.
\end{equation}
\end{theorem}

It was observed earlier in \cite{Sball} that the special case of \eqref{5.5} with $r=\la^{-1}$ follows from results of
B\'erard~\cite{PB}.  Of course it would be very interesting to see if one could establish a result like \eqref{5.4} for $r$ all the way
up to a logarithmic scale.  The reason that we can obtain \eqref{5.5} for the range of $r$ there is that the radii are sufficiently small
so that, presumably, only functions with the profiles of zonal functions or highest weight spherical harmonics can saturate 
\eqref{5.4} for such $r$, and  eigenfunctions with their profiles cannot exist on manifolds of nonpositive curvature.

The proof of Theorem~\ref{ball} is based on a simple variation on arguments in \cite{BS2}, \cite{BSlog} and \cite{SZKN}.

First, to prove \eqref{5.5}, we note that if $\rho\in {\mathcal S}(\R)$ satisfies $\rho(0)=0$ and if $P=\sqrt{-\Delta_g}$, then
for any $T\ge1$
\begin{equation}\label{reproduce}\rho(T(\la-P))e_\la = e_\la.\end{equation}
We shall also want to assume, as we may, that the Fourier transform of $\rho$ satisfies $\Hat \rho(t)=0$ if $|t|\ge 1/2$ and we shall
take $T=T(\la)=c\log\la$, where $c$ depends on $(M,g)$ and $\e$ and will be specified in a moment.

Based on \eqref{reproduce}, we see that we would obtain \eqref{5.5} if we could show that 
\begin{equation}\label{5.4'}\tag{5.4$'$}
\bigl\|\rho(T(\la-P))f\bigr\|_{L^2(B(x,r))}\le C_\e \bigl(r/\log\la\bigr)^{\frac12}\|f\|_{L^2(M)}, \quad \la^{-1}\le r\le \la^{-\frac12-\e}.
\end{equation}
By a routine $TT^*$ argument, \eqref{5.4'} is equivalent to showing that
\begin{multline}\label{5.4''}\tag{5.4$''$}
\bigl\|\chi(T(\la-P))h\|_{L^2(B(x,r))}\le C_\e \bigl(r/\log \la\bigr)\|h\|_{L^2},  
\\
 \text{if } \, \text{supp }h\subset B(x,r), \, \,
\text{and } \, \, \, \la^{-1}\le r\le \la^{-\frac12-\e},
\end{multline} 
with
$\chi=|\rho|^2$.  Note that because of our assumptions for $\rho$ we have
\begin{equation}\label{5.6}
\Hat \chi(t)=0 \quad \text{if } \, \, |t|\ge 1.
\end{equation}

Notice also that
\begin{align*}
\chi(T(\la-P))&=\frac1{2\pi T}\int \Hat \chi(t/T) e^{it\la} e^{-itP} \, dt
\\
&=\frac1{\pi T}\int \Hat \chi(t/T) e^{i\la t} \cos t\sqrt{-\Delta_g} \, dt \, +\, \chi(T(\la+P)).
\end{align*}
Since $\chi\in {\mathcal S}(\R)$ and $P$ is nonnegative, it follows that $\chi(T(\la+P))$ has a kernel which
is $O(\la^{-N})$ for all $N$.  Therefore, we would have \eqref{5.4''} if we could show that for $T=c\log\la$
(and appropriate $c>0$) we have
\begin{multline}\label{5.7}
T^{-1}\Bigl\| \int \Hat \chi(t/T) e^{i\la t} \cos t\sqrt{-\Delta_g} \, h\, dt\Bigr\|_{L^2(B(x,r))}\le C_\e \bigl(r/\log\la\bigr)\|h\|_{L^2},
\\
\text{if } \, \, \text{supp }h\subset B(x,r), \, \, \, \la^{-1}\le r\le \la^{-\frac12-\e}.
\end{multline}

If $(\Rn,\tilde g)$ denotes the universal cover of $(M,g)$  and if $D\subset \Rn$ is an associated Dirichlet domain then, as in 
\cite{BS2}, \cite{BSlog} and \cite{SZKN} we shall use the fact that
$$\bigl(\cos t\sqrt{-\Delta_g}\bigr)(x,y) = \sum_{\alpha\in \Gamma}\bigl(\cos t\sqrt{-\Delta_{\tidle g}}\bigr)(\tilde x,\alpha(\tilde y)),$$
where $\Gamma$ are the deck transformations coming from our covering and $\tilde x,\tilde y\in D$ denote the
lifts of our $x,y\in M$, respectively, to the fundamental domain.

Thus, the operator in \eqref{5.7} can be rewritten as $\sum_{\alpha\in \Gamma}\chi_\alpha$ where the kernel of $\chi_\alpha$ is
$$\chi_\alpha(x,y)=T^{-1}\int \Hat \chi(t/T) e^{i\la t} \bigl(\cos t\sqrt{-\Delta_{\tilde g}}\bigr)(\tilde x,\alpha(\tilde y)) \, dt.$$
By (3.8) in \cite{BSlog}
$$|\chi_\alpha(x,y)|\le CT^{-1}\la^{\frac{n-1}2}, \quad \text{if } \, \, \alpha \ne {\rm Identity},$$
and so, by Young's inequality,
\begin{equation}\label{5.8}
\|\chi_\alpha h\|_{L^2(B(x,r))}\le CT^{-1}\la^{\frac{n-1}2}r^n\|h\|_{L^2}, \quad \text{if } \, \, \text{supp }h\subset B(x,r).
\end{equation}
Note that 
$$\la^{\frac{n-1}2}r^n\le r\la^{-(n-1)\e}, \quad \text{if } \, \, r\le \la^{-\frac12-\e}.$$
Also, by Huygen's principle,  \eqref{5.6} and volume comparison theorems there are $O(\exp(C_MT))$ nonzero terms
$\chi_\alpha$ for some constant $C_M$ depending on $(M,g)$.  Therefore, 
since we are taking $T=c\log\la$, we conclude that 
 if we choose $c>0$ 
to be small enough then
we can sum the bounds in \eqref{5.8} to 
obtain
\begin{equation}\label{5.8'}\tag{5.8$'$}
\sum_{{\rm Identity}\ne \alpha\in \Gamma}\|\chi_\alpha\|_{L^2(B(x,r))\to L^2(B(x,r))} \le \la^{-\frac{\e}2}r,
\end{equation}
which is better than the bounds posited in \eqref{5.7}.
As a result, we would obtain this inequality if we also knew that
$$\|\chi_\alpha\|_{L^2(B(x,r))\to L^2(B(x,r))}\le T^{-1}r, \quad \text{if} \, \, \, \alpha={\rm Identity}.$$
But this follows from routine arguments.  Indeed, if follows from a straightforward modification of the argument
in \cite{Sball} that was used to establish \eqref{5.4}.

\section*{Acknowledgements}
We thank Matthew Blair and Steve Zelditch for helpful comments and suggestions.

\end{document}